\documentclass{amsart}

\usepackage{amsmath,amsthm,amsfonts,amscd,amssymb}
\usepackage{verbatim}
\usepackage{color} 

\newtheorem{thm}{Theorem}[section]
\newtheorem{cor}[thm]{Corollary}
\newtheorem{lem}[thm]{Lemma}
\newtheorem{prop}[thm]{Proposition}

\theoremstyle{definition}

\theoremstyle{remark}
\newtheorem{rem}{Remark}[section]

\begin{document}

\title{Height bounds on zeros of quadratic forms over~$\overline{\mathbb Q}$}

\author{ Lenny Fukshansky}\thanks{The author was partially supported by the Simons Foundation grant \#279155 and NSA grant \#130907.}



\address{Department of Mathematics, 850 Columbia Avenue, Claremont McKenna College, Claremont, CA 91711}
\email{lenny@cmc.edu}



\begin{abstract}
In this paper we establish three results on small-height zeros of quadratic polynomials over~$\overline{\mathbb Q}$. For a single quadratic form in $N \geq 2$ variables on a subspace of $\overline{\mathbb Q}^N$, we prove an upper bound on the height of a smallest nontrivial zero outside of an algebraic set under the assumption that such a zero exists. For a system of $k$ quadratic forms on an $L$-dimensional subspace of~$\overline{\mathbb Q}^N$, $N \geq L \geq \frac{k(k+1)}{2}+1$, we prove existence of a nontrivial simultaneous small-height zero. For a system of one or two inhomogeneous quadratic and $m$ linear polynomials in $N \geq m+4$ variables, we obtain upper bounds on the height of a smallest simultaneous zero, if such a zero exists. Our investigation extends previous results on small zeros of quadratic forms, including Cassels' theorem and its various generalizations and contributes to the literature of so-called ``absolute" Diophantine results with respect to height. All bounds on height are explicit.

\end{abstract}

\maketitle

\def\A{{\mathcal A}}
\def\AA{{\mathfrak A}}
\def\B{{\mathcal B}}
\def\C{{\mathcal C}}
\def\D{{\mathcal D}}
\def\E{{\mathcal E}}
\def\F{{\mathcal F}}
\def\Ff{{\mathfrak F}}
\def\G{{\mathcal G}}
\def\x{{\mathcal H}}
\def\I{{\mathcal I}}
\def\J{{\mathcal J}}
\def\K{{\mathcal K}}
\def\kk{{\mathfrak K}}
\def\L{{\mathcal L}}
\def\LL{{\mathfrak L}}
\def\M{{\mathcal M}}
\def\O{{\mathcal O}}
\def\W{{\omega}}
\def\CC{{\mathfrak C}}
\def\mm{{\mathfrak m}}
\def\MM{{\mathfrak M}}
\def\OO{{\mathfrak O}}
\def\P{{\mathcal P}}
\def\R{{\mathcal R}}
\def\s{{\mathcal S}}
\def\V{{\mathcal V}}
\def\X{{\mathcal X}}
\def\XX{{\mathfrak X}}
\def\Y{{\mathcal Y}}
\def\Z{{\mathcal Z}}
\def\H{{\mathcal H}}
\def\cee{{\mathbb C}}
\def\pee{{\mathbb P}}
\def\que{{\mathbb Q}}
\def\real{{\mathbb R}}
\def\zed{{\mathbb Z}}
\def\hyp{{\mathbb H}}
\def\aaa{{\mathbb A}}
\def\ff{{\mathbb F}}
\def\kk{{\mathfrak K}}
\def\qbar{{\overline{\mathbb Q}}}
\def\kbar{{\overline{K}}}
\def\ybar{{\overline{Y}}}
\def\kkbar{{\overline{\mathfrak K}}}
\def\ubar{{\overline{U}}}
\def\eps{{\varepsilon}}
\def\ahat{{\hat \alpha}}
\def\bhat{{\hat \beta}}
\def\gt{{\tilde \gamma}}
\def\h{{\tfrac12}}
\def\dd{{\partial}}
\def\baa{{\boldsymbol \alpha}}
\def\bfa{{\boldsymbol a}}
\def\bfb{{\boldsymbol b}}
\def\be{{\boldsymbol e}}
\def\bei{{\boldsymbol e_i}}
\def\bff{{\boldsymbol f}}
\def\bc{{\boldsymbol c}}
\def\bm{{\boldsymbol m}}
\def\bk{{\boldsymbol k}}
\def\bi{{\boldsymbol i}}
\def\bl{{\boldsymbol l}}
\def\bq{{\boldsymbol q}}
\def\bu{{\boldsymbol u}}
\def\bt{{\boldsymbol t}}
\def\bs{{\boldsymbol s}}
\def\bfu{{\boldsymbol u}}
\def\bv{{\boldsymbol v}}
\def\bw{{\boldsymbol w}}
\def\bx{{\boldsymbol x}}
\def\bX{{\boldsymbol X}}
\def\bz{{\boldsymbol z}}
\def\bwy{{\boldsymbol y}}
\def\bY{{\boldsymbol Y}}
\def\bL{{\boldsymbol L}}
\def\ba{{\boldsymbol a}}
\def\bb{{\boldsymbol b}}
\def\bet{{\boldsymbol\eta}}
\def\bxi{{\boldsymbol\xi}}
\def\bo{{\boldsymbol 0}}
\def\bol{{\boldkey 1}_L}
\def\ep{\varepsilon}
\def\p{\boldsymbol\varphi}
\def\q{\boldsymbol\psi}
\def\rank{\operatorname{rank}}
\def\aut{\operatorname{Aut}}
\def\lcm{\operatorname{lcm}}
\def\sgn{\operatorname{sgn}}
\def\spn{\operatorname{span}}
\def\md{\operatorname{mod}}
\def\Norm{\operatorname{Norm}}
\def\dim{\operatorname{dim}}
\def\det{\operatorname{det}}
\def\Vol{\operatorname{Vol}}
\def\rk{\operatorname{rk}}
\def\ord{\operatorname{ord}}
\def\ker{\operatorname{ker}}
\def\div{\operatorname{div}}
\def\Gal{\operatorname{Gal}}
\def\GL{\operatorname{GL}}
\def\p{\operatorname{p}}
\def\q{\operatorname{q}}
\def\t{\operatorname{t}}
\def\hs{{\hat \sigma}}
\def\chr{\operatorname{char}}

\section{Introduction and statement of results}

Results over~$\qbar$, the algebraic closure of~$\que$, have received some attention in Diophantine geometry and theory of height functions, especially in the recent years. A classical example of such a result is the ``absolute" version of Siegel's lemma of Roy and Thunder~\cite{absolute:siegel}, which can be viewed as a statement about existence of small-height solutions to systems of homogeneous linear equations over~$\qbar$. The Roy-Thunder result extends the Bombieri-Vaaler version of Siegel's lemma~\cite{vaaler:siegel} over number fields to the $\qbar$ situation, producing height bounds independent of any number field, hence the name absolute.

To move beyond the linear equations, the investigation of small-height zeros of quadratic forms was initiated in the celebrated paper of Cassels \cite{cassels:small}, and later continued by a number of authors (see \cite{cassels_overview} for a detailed overview). Most of the work here has been done over fixed number fields and function fields, however some absolute results have also been produced. For instance, Vaaler's theorem~\cite{vaaler:smallzeros} on small-height maximal totally isotropic subspaces over a fixed number field has been extended over~$\qbar$ in~\cite{quad_qbar}. Techniques used to produce absolute results often differ from the methods over a fixed number field, since sets of points of bounded height are no longer finite (i.e., the Northcott property fails). In a recent paper~\cite{cfh}, we have obtained results on existence of zeros and isotropic subspaces of a quadratic space outside of a finite union of varieties over a fixed number field or global function field. The first goal of the present paper is to obtain analogous results over~$\qbar$.

Let $J \geq 1$ be an integer. For each $1 \leq i \leq J$, let $k_i \geq 1$ be an integer and let
$$P_{i1}(X_1,\dots,X_N),\dots,P_{ik_i}(X_1,\dots,X_N)$$
be homogeneous polynomials of respective degrees $m_{i1},\dots,m_{ik_i} \geq 1$. Let
$$Z(P_{i1},\dots,P_{ik_i}) = \{ \bx \in \qbar^N : P_{i1}(\bx) = \dots = P_{ik_i}(\bx) = 0 \},$$
and define the algebraic set
\begin{equation}
\label{Z_K}
\Z = \bigcup_{i=1}^J Z(P_{i1},\dots,P_{ik_i}).
\end{equation}
For each $1 \leq i \leq J$ let $M_i = \max_{1 \leq j \leq k_i} m_{ij}$, and define
\begin{equation}
\label{def_M}
M = M(\Z) := \sum_{i=1}^J M_i.
\end{equation}
The basic notation of the arithmetic theory of quadratic forms, which is used in the statement of Theorem~\ref{miss_hyper}, is reviewed below in Section~\ref{notation}, along with definitions of the appropriate height functions.
\smallskip

\begin{thm} \label{miss_hyper} Let $F$ be a quadratic form in $N$ variables over $\qbar$, $V$ be an $L$-dimensional subspace of $\qbar^N$, and $m$ be the dimension of a maximal totally isotropic subspace of the quadratic space $(V,F)$.  Suppose that $F$ has a nontrivial zero in $V\setminus \Z$.  Then there exist $m$ linearly independent zeros $\bx_1,\dots,\bx_m$ of $F$ in $V\setminus \Z$ such that
\begin{equation}
\label{h_order}
H(\bx_1) \leq H(\bx_2) \leq \dots \leq H(\bx_m),\ h(\bx_1) \leq h(\bx_2) \leq \dots \leq h(\bx_m),
\end{equation}
and for each $1 \leq n \leq m$,
\begin{equation}
\label{miss_hyper_bnd}
H(\bx_n) \leq h(\bx_n) \leq T(L,M+1) H(F)^{\max \{r,29/2\}} \H(V)^{30},
\end{equation}
where $r$ is the rank of $F$ on $V$ and the constant $T(L,M)$ is defined by \eqref{T_const_qb} below.
\end{thm}

As a corollary of Theorem~\ref{miss_hyper}, we can also obtain a statement on the existence of nested sequences of totally isotropic subspaces of $(V,F)$ of bounded height not contained in the algebraic set~$\Z$.

\begin{cor} \label{miss_subspace} Let $r$ and $\W = [r/2]$ be the rank and the Witt index, respectively, of the quadratic space $(V,F)$ in Theorem \ref{miss_hyper}.  Then for each pair of indices $(n,k)$ with $1 \leq n,k \leq m$ there exists a totally isotropic subspace $W^k_n$ of $(V,F)$ such that $\dim_K W^k_n = k$,
$ W^1_n \subset W^2_n \subset \dots \subset W^m_n$, and so $W^k_n \nsubseteq \Z$ for each $1 \leq k \leq m$; also
\begin{equation}
\label{cor_bnd}
\H(W^m_n) \leq T^1(L,M,N,\W) H(F)^{\t_1(\W,r)} \H(V)^{\t_2(\W)},
\end{equation}
where the exponents are given by
\begin{eqnarray}
\label{t_exp}
& & \t_1(\W,r) = (\W-1)^2 + \frac{4\W + 4 + (4\W+7)\max \{r,29/2\}}{3} + r,\nonumber \\
& & \t_2(\W) = 40\W + 70 + \frac{4\W+4}{3},
\end{eqnarray}
and the constant $T^1(L,M,N,m)$ is defined in~\eqref{cor_bnd_const_qb} below. In addition, for each $1 \leq k < m$,
\begin{equation}
\label{cor_bnd_1}
\H(W^k_n) \leq 3^{\frac{m(m-1)}{4}} N^{\frac{k}{2}} H(\bx_n) \H(W^m_n),
\end{equation}
where $\bx_n$ is from Theorem \ref{miss_hyper}.
\end{cor}
\smallskip

\noindent
The main line of our argument in the proof of Theorem~\ref{miss_hyper} is analogous to that of~\cite{cfh}, however one of the tools we need in the proof of Theorem~\ref{miss_hyper} and Corollary~\ref{miss_subspace} is the result of~\cite{quad_qbar} on the existence of a small-height maximal totally isotropic subspace of a quadratic space over~$\qbar$. 

While working in the ``absolute" setting presents some difficulties (like the failure of Northcott property), one naturally expects some problems to become easier over~$\qbar$. For instance, a version of the original Cassels' theorem can be proved in a much simpler way with a considerably better bound over~$\qbar$. Specifically, here is Lemma~4.1 of~\cite{quad_qbar}.

\begin{lem} [Lemma 4.1, \cite{quad_qbar}] \label{zero} Let $2 \leq L \leq N$ and $V \subseteq \qbar^N$ an $L$-dimensional subspace. Let $F(\bX)$ be a quadratic form in $N$ variables over $\qbar$. There exists $\bo \neq \bwy \in V$ such that $F(\bwy)=0$ and
\begin{equation}
\label{zero_eq}
h(\bwy) \leq 8 \times 3^{2(L-1)} \H(V)^{\frac{4}{L}} H(F)^{\frac{1}{2}}.
\end{equation}
\end{lem}

\noindent
It is not known if a bound for the height of a simultaneous zero of a system of quadratic forms on a subspace $V \subseteq K^N$ exists over a fixed number field~$K$: this question is connected with a very general version of Hilbert's 10th problem over number fields, and currently appears to be out of reach. On the other hand, over $\qbar$ the problem is more tractable, as we show next.

\begin{thm} \label{k_forms} Let $k \geq 2$ be an integer, $F_1,\dots,F_k$ be quadratic forms in $N$ variables over $\qbar$, and let $V \subseteq \qbar^N$ be an $L$-dimensional subspace, $N \geq L \geq \frac{k(k+1)}{2} + 1$. There exists $\bo \neq \bw \in V$ such that $F_m(\bw) = 0$ for all $1 \leq m \leq k$ and
\begin{equation}
\label{hw_main}
h(\bw) \leq \left( 3^{\frac{L^2}{2}} N^{\frac{3(L+1)}{2}} \H(V) \right)^{20B_k^2/81} \left( \prod_{n=1}^{k-1} H(F_n) \right)^{B_k} H(F_k)^2,
\end{equation}
where $B_2=9$ and
\begin{equation}
\label{BK}
B_k = \frac{1}{4} \times 36^{2^{k-2}} \prod_{m=3}^k m^{2^{k-m+1}}
\end{equation}
for all $k \geq 3$.
\end{thm}

\noindent
Notice that the exponents on heights of $V$ and the quadratic forms in Theorem~\ref{k_forms} depend only on $k$, the number of forms, not on their number of variables or dimension of the space, as they usually do in bounds over global fields.

Finally, we obtain a bound on the height of a nontrivial solution of a system of one or two quadratic and a collection of linear inhomogeneous equations over~$\qbar$.

\begin{thm} \label{thm_inh} Let $F$ and $G$ be quadratic polynomials in $N \geq 4$ variables over $\qbar$, possibly inhomogeneous. Let $m$ be an integer, $0 \leq m \leq N-4$, and $\L_1,\dots,\L_m$ be linear polynomials in $N$ variables over $\qbar$, possibly inhomogeneous; the case $m=0$ just means that there are no linear polynomials. Suppose that the system
\begin{equation}
\label{system-2}
F(\bx) = G(\bx) = \L_1(\bx) = \dots = \L_m(\bx) = 0
\end{equation}
has a nontrivial solution in~$\qbar^N$. Then there exists a point $\bo \neq \bwy \in \qbar^N$ such that
$$F(\bwy) = \L_1(\bwy) = \dots = \L_m(\bwy) = 0$$
and
\begin{equation}
\label{inh_bnd_1}
h(\bwy) \leq 8 (N+1)^{2m} 3^{2(N-m+1)(N-m)} H(F)^{\frac{1}{2}} \prod_{i=1}^m H(\L_i)^4.
\end{equation}
There also exists a point $\bo \neq \bw \in \qbar^N$ such that
$$F(\bw) = G(\bw) = \L_1(\bw) = \dots = \L_m(\bw) = 0$$
and
\begin{equation}
\label{inh_bnd_2}
h(\bw) \leq \M(m,N) H(F)^{58} H(G)^3 \prod_{i=1}^m H(\L_i)^{180},
\end{equation}
where
\begin{equation}
\label{MmN-const}
\M(m,N) = 18 \times 8^{38} (N+1)^{90m+8} (N+1-m)^{36} 3^{90(N-m+1)(N-m)}.
\end{equation}
\end{thm}
\smallskip

\begin{rem} \label{over_K} In fact, our method can also be used to obtain analogues of Theorems~\ref{k_forms} and~\ref{thm_inh} with points in question having bounded degree over a fixed number field. By the Northcott property, this provides actual search bounds for zeros of systems of quadratic and linear equations as above. On the other hand, the bounds on height we can obtain this way are weaker.
\end{rem}

This paper is organized as follows. In Section~\ref{notation} we set the notation, define the height functions, and review the basic terminology in the theory of quadratic forms. We prove Theorem~\ref{miss_hyper} along with Corollary~\ref{miss_subspace} in Section~\ref{section_miss_hyper}. In Section~\ref{section_k_forms} we prove Theorem~\ref{k_forms}. Finally, Theorem~\ref{thm_inh} is proved in Section~\ref{inh}.
\bigskip

\section{Preliminaries} \label{notation}

\subsection{Notations and heights}
We start with some notation, following \cite{null}, \cite{witt}, and~\cite{quad_qbar} (see also~\cite{bombieri_gubler} for a comprehensive overview of the theory of height functions and their many properties). Let $K$ be a number field and let $d = [K:\que]$ be the global degree of $K$ over $\que$. Let $M(K)$ be the set of all places of $K$. For each place $v \in M(K)$ we write $K_v$ for the completion of $K$ at $v$ and let $d_v = [K_v:\que_v]$ be the local degree of $K$ at $v$.  For each place $v \in M(K)$ we define the absolute value $|\ |_v$ to be the unique absolute value on $K_v$ that extends either the usual absolute value on $\real$ or $\cee$ if $v | \infty$, or the usual $p$-adic absolute value on $\que_p$ if $v|p$, where $p$ is a prime. With this choice of absolute values, the product formula reads as follows:
\begin{equation}
\label{product_formula}
\prod_{v \in M(K)} |a|^{d_v}_v = 1, \quad \mbox{ for all $a \in K^\times$}.
\end{equation}

For each $v \in M(K)$ define a local height $H_v$ on $K_v^N$ by
$$H_v(\bx) = \max_{1 \leq i \leq N} |x_i|^{d_v}_v,$$
for each $\bx \in K_v^N$. For each $v | \infty$ we also define another local height
$$\H_v(\bx) = \left( \sum_{i=1}^N |x_i|_v^2 \right)^{d_v/2}.$$
Notice that these local heights $H_v$ and $\H_v$ are $\ell^{\infty}$- and $\ell^2$-norms, respectively, on the vector spaces~$K_v^N$. Then we can define two slightly different global height functions on $K^N$:
\begin{equation}
\label{global_heights}
H(\bx) = \left( \prod_{v \in M(K)} H_v(\bx) \right)^{1/d},\ \ \H(\bx) = \left( \prod_{v \nmid \infty} H_v(\bx) \times \prod_{v | \infty} \H_v(\bx) \right)^{1/d},
\end{equation}
for each $\bx \in K^N$. These height functions are {\it homogeneous}, in the sense that they are defined on the projective space over $K^N$ thanks to the product formula (\ref{product_formula}).  It is easy to see that
\begin{equation}
\label{ht_ineq_sqrt}
H(\bx) \leq \H(\bx) \leq \sqrt{N} H(\bx).
\end{equation}
We also define the {\it inhomogeneous} height
$$h(\bx) = H(1,\bx),$$
which generalizes the Weil height on algebraic numbers.  Clearly, $h(\bx) \geq H(\bx)$ for each $\bx \in K^N$.  We extend the height functions $H_v$, $\H_v$, $H$ and $h$ to polynomials by evaluating the height of their coefficient vectors, and to matrices by viewing them as vectors.  However, if $X$ is a matrix with $\bx_1, \ldots, \bx_L$ as its columns, then $\H(X)$ will always denote the height $\H(\bx_1\wedge \dots \wedge \bx_L)$, where $\bx_1 \wedge \dots \wedge \bx_L$ is viewed as a vector in $K^{\binom{N}{L}}$ under the standard embedding.

Let $V$ be an $L$-dimensional subspace of $K^N$.  Then there exist $N \times L$ matrix $X$ and $(N-L) \times N$ matrix $A$, both are over $K$, such that
$$V = \{ X \bt : \bt \in K^L \} = \{ \bx \in K^N : A \bx = 0 \}.$$
The Brill-Gordan duality principle \cite{gordan:1} (also see Theorem 1 on p. 294 of \cite{hodge:pedoe}) implies that $\H(X)=\H(A^t)$, and $\H(V)$ is defined to be this common value.  This coincides with the choice of heights in \cite{vaaler:siegel}.

An important observation is that due to the normalizing exponent $1/d$ in (\ref{global_heights}) all our heights are {\it absolute}, meaning that they do not depend on the number field of definition.  Therefore we have defined the necessary height functions over~$\qbar$. We also recall here the {\it Northcott property}, satisfied by our height functions. Given $\bo \neq \bx = (x_1,\dots,x_N) \in \qbar^N$, let us write $[\bx]$ for the corresponding projective point. For a fixed number field $K$, define $\deg_K \bx := [K(x_1,\dots,x_N):K]$ and $\deg_K [\bx] := \min \{ \deg_K \bwy : \bwy \in \qbar^N, [\bwy] = [\bx] \}$.
Let $C \in \real_{>0}$, $D \in \zed_{>0}$, $K$ a fixed number field. Then the sets
$$\left\{ \bx \in \qbar^N : \deg_K \bx \leq D, h(\bx) \leq C \right\},\ \left\{ [\bx] \in \pee(\qbar^N) : \deg_K [\bx] \leq D, H(\bx) \leq C \right\}$$
are finite.

We will also need a few technical lemmas detailing some basic properties of heights. The first one bounds the height of a linear combination of vectors (see, for instance, Lemma~2.1 of~\cite{cfh-1}).

\begin{lem} \label{sum_height} For $\xi_1,...,\xi_L \in K$ and $\bx_1,...,\bx_L \in K^N$,
$$H \left( \sum_{i=1}^L \xi_i \bx_i \right) \leq h \left( \sum_{i=1}^L \xi_i \bx_i \right) \leq L h(\bxi) \prod_{i=1}^L h(\bx_i),$$
where $\bxi = (\xi_1,...,\xi_L) \in K^L$.
\end{lem}

\noindent
The second one is an adaptation of Lemma 4.7 of \cite{absolute:siegel} to our choice of height functions, using~\eqref{ht_ineq_sqrt}.

\begin{lem} \label{lem_4.7} Let $V$ be a subspace of $K^N$, $N \geq 2$, and let subspaces $U_1,\dots,U_n \subseteq V$ and vectors $\bx_1,\dots,\bx_m \in V$ be such that
$$V = \spn_K \{ U_1,\dots,U_n,\bx_1,\dots,\bx_m \}.$$
Then
$$\H(V) \leq N^{m/2} \H(U_1) \dots \H(U_n) H(\bx_1) \dots H(\bx_m).$$
\end{lem}

\noindent
The next one is an adaptation of Lemma 2.3 of \cite{witt} to our choice of height functions, using~\eqref{ht_ineq_sqrt}.

\begin{lem} \label{lem_2.3} Let $X$ be a $J \times N$ matrix over $K$ with row vectors $\bx_1,...,\bx_J$, and let $F$ be a symmetric bilinear form in $2N$ variables over $K$ (we also write $F$ for its $N \times N$ coefficient matrix). Then
$$\H(X F) \leq N^{3J/2} H(F)^J \prod_{i=1}^J H(\bx_i).$$
\end{lem}

\noindent
The next one is a bound on the absolute values of a bilinear form at a given pair of vectors.

\begin{lem} \label{F_xy} Let $\bx,\bwy \in K^N$ and let $F$ be a symmetric bilinear form in $2N$ variables over $K$ (we also write $F = (f_{ij})_{1 \leq i,j \leq N}$ for its $N \times N$ coefficient matrix). For each $v \in M(K)$, we have:
$$|F(\bx,\bwy)|^{d_v}_v \leq \left\{ \begin{array}{ll}
N^{2d_v} H_v(F) H_v(\bx) H_v(\bwy) & \mbox{if $v \mid \infty$,} \\
H_v(F) H_v(\bx) H_v(\bwy) & \mbox{if $v \nmid \infty$,}
\end{array}
\right.$$
and hence the Weil height of $F(\bx,\bwy)$ is
$$h(F(\bx,\bwy)) = \prod_{v \in M(K)} \max \left\{ 1, |F(\bx,\bwy)|^{d_v}_v \right\}^{1/d} \leq N^2 H(F) H(\bx) H(\bwy).$$
\end{lem}

\proof
Notice that
$$F(\bx,\bwy) = \sum_{i=1}^N \sum_{j=1}^N f_{ij} x_iy_j.$$
By triangle inequality, for each $v \mid \infty$,
$$|F(\bx,\bwy)|_v \leq \sum_{i=1}^N \sum_{j=1}^N |f_{ij} x_iy_j|_v \leq N^2 \max_{1 \leq i,j \leq N} |f_{ij}|_v \max_{1 \leq i \leq N} |x_i|_v \max_{1 \leq j \leq N} |y_j|_v.$$
On the other, for $v \nmid \infty$ the ultrametric inequality implies
$$|F(\bx,\bwy)|_v \leq \max_{1 \leq i,j \leq N} |f_{ij} x_iy_j|_v \leq \max_{1 \leq i,j \leq N} |f_{ij}|_v \max_{1 \leq i \leq N} |x_i|_v \max_{1 \leq j \leq N} |y_j|_v.$$
\endproof

\noindent
The next one is Lemma 2.2 of \cite{witt}.

\begin{lem} \label{intersection} Let $U_1$ and $U_2$ be subspaces of $K^N$. Then
\begin{equation}
\label{inter_eq}
\H(U_1 \cap U_2) \leq \H(U_1) \H(U_2).
\end{equation}
\end{lem}
\smallskip
 
\begin{rem} \label{strong_intersection} It should be remarked that a stronger version of inequality~\eqref{inter_eq} has been produced in~\cite{schmidt-remark} and~\cite{vaaler:struppeck}, specifically:
\begin{equation}
\label{inter_eq-strong}
\H(U_1+U_2) \H(U_1 \cap U_2) \leq \H(U_1) \H(U_2).
\end{equation}
Unfortunately, for our purposes using the stronger inequality~\eqref{inter_eq-strong} instead of~\eqref{inter_eq} does not seem to have an immediate benefit: while the bounds come out considerably more complicated and hard to read, it is not clear how much better they are, since the quantity $\H(U_1+U_2)$ is usually hard to nontrivially estimate from below.
\end{rem}

\begin{rem} \label{overkbar} Lemmas~\ref{sum_height} - \ref{intersection} also hold verbatim with $K$ replaced by $\qbar$.
\end{rem}
\smallskip

\subsection{Quadratic Forms}

Here we introduce some basic language of quadratic forms which are necessary for subsequent discussion.  For an introduction to the subject, the readers are referred to, for instance, Chapter 1 of \cite{scharlau}.  For the sake of more generality, we allow $K$ to be any field of characteristic not 2.  We write
$$F(\bX,\bY) = \sum_{i=1}^N \sum_{j=1}^N f_{ij} X_i Y_j$$
for a symmetric bilinear form in $2N$ variables with coefficients $f_{ij} = f_{ji}$ in $K$, and $F(\bX) = F(\bX,\bX)$ for the associated quadratic form in $N$ variables; we also use $F$ to denote the symmetric $N \times N$ coefficient matrix $(f_{ij})_{1 \leq i,j \leq N}$.    Let $V$ be an $L$-dimensional subspace of $K^N$, $2 \leq L \leq N$.  Then $F$ is also defined on $V$, and we write $(V,F)$ for the corresponding quadratic space.

A point $\bx$ in a subspace $U$ of $V$ is called singular if $F(\bx,\bwy) = 0$ for all $\bwy \in U$, and it is called nonsingular otherwise. For each subspace $U$ of $(V,F)$, its radical is the set
$$U^{\perp} := \{ \bx \in U : F(\bx, \bwy) = 0\ \forall\ \bwy \in U \}$$
which is the subspace of all singular points in $U$.   We define $\lambda(U):=\dim_K U^{\perp}$, and will write $\lambda$ to denote $\lambda(V)$. A subspace $U$ of $(V,F)$ is called regular if $\lambda(U)=0$.

A point $\bo \neq \bx \in V$ is called isotropic if $F(\bx)=0$ and anisotropic otherwise. A subspace $U$ of $V$ is called isotropic if it contains an isotropic point, and it is called anisotropic otherwise. A totally isotropic subspace $W$ of $(V,F)$ is a subspace such that for all $\bx,\bwy \in W$, $F(\bx,\bwy)=0$. All maximal totally isotropic subspaces of $(V,F)$ contain $V^{\perp}$ and have the same dimension. Given any maximal totally isotropic subspace $W$ of $V$, let
$$\W = \W(V) := \dim_K(W)-\lambda,$$
which is the Witt index of $(V,F)$.  If $K=\kbar$, then $\W= [(L-\lambda)/2]$, where $[\ ]$ stands for the integer part function.

If two subspaces $U_1$ and $U_2$ of $(V,F)$ are orthogonal, we write $U_1 \perp U_2$ for their orthogonal sum. If $U$ is a regular subspace of $(V,F)$, then $V = U \perp \left( \perp_V(U) \right)$ and $U \cap  \left( \perp_V(U) \right) = \{\boldsymbol 0\}$, where
\begin{equation}
\label{perp_V}
\perp_V(U) := \{ \bx \in V : F(\bx,\bwy) = 0\ \forall\ \bwy \in U \}
\end{equation}
is the orthogonal complement of $U$ in $V$. Two vectors $\bx,\bwy \in V$ are called a hyperbolic pair if $F(\bx) = F(\bwy) = 0$ and $F(\bx,\bwy) \neq 0$; the subspace $\hyp(\bx,\bwy) := \spn_K \{\bx,\bwy\}$ that they generate is regular and is called a hyperbolic plane. An orthogonal sum of hyperbolic planes is called a hyperbolic space. Every hyperbolic space is regular. It is well known that there exists an orthogonal Witt decomposition of the quadratic space $(V,F)$ of the form
\begin{equation}
\label{decompose}
V = V^{\perp} \perp \hyp_1 \perp\ \dots \perp \hyp_{\W} \perp U,
\end{equation}
where $\hyp_1, \dots, \hyp_{\W}$ are hyperbolic planes and $U$ is an anisotropic subspace, which is determined uniquely up to isometry. The rank of $F$ on $V$ is $r:=L-\lambda$. In case $K=\kbar$, $\dim_K U = 1$ if $r$ is odd and 0 if $r$ is even. Therefore a regular even-dimensional quadratic space over $\qbar$ is always hyperbolic.

We are now ready to proceed.
\bigskip

\section{Isotropic points missing varieties}
\label{section_miss_hyper}

In this section we prove Theorem \ref{miss_hyper} and Corollary~\ref{miss_subspace}. We start with a non-vanishing lemma for polynomials, which is the direct analogue of Lemma~4.1 of~\cite{cfh}.

\begin{lem} \label{nonvanish} Let $N, M \geq 1$ be integers and let $P(\bX):= P(X_1,\dots,X_N) \in \qbar[X_1,\dots,X_N]$ be a polynomial that is not identically zero with $\deg P \leq M$. Then there exists $\bz \in \qbar^N$ such that $P(\bz) \neq 0$ and
$$h(\bz) \leq 1.$$
\end{lem}

\proof
The conclusion of the lemma follows immediately from Lemma 4.1 of \cite{null} combined with the argument in Section~7 of \cite{null} (in particular, see formulas (44) and (45) of \cite{null}).
\endproof

We also need a technical lemma providing a bound on the height of a restriction of a polynomial to a subspace.

\begin{lem} \label{restrict} Let $N, M \geq 1$ be integers and let $P(\bX) := P(X_1,\dots,X_N) \in \qbar[X_1,\dots,X_N]$ be a polynomial of degree $M$. Let $V \subseteq \qbar^N$ be an $L$-dimensional subspace, $1 \leq L \leq N$, such that $P$ is not identically zero on $V$. Let $\bx_1,\dots,\bx_L$ be a basis for $V$ over $\qbar$, write $A$ for the $N \times L$ basis matrix $(\bx_1,\dots,\bx_L)$, and define
$$P_A(Y_1,\dots,Y_L) = P(Y_1\bx_1 + \dots + Y_L\bx_L) \in \qbar[Y_1,\dots,Y_L],$$
so that $P_A$ is a restriction of $P$ to $V$. Then $P_A$ is a polynomial of degree $M$ in $L$ variables over $\qbar$, and
\begin{equation}
\label{ht_PA}
H(P_A) \leq L^{M} H(P) \prod_{i=1}^L h(\bx_i)^M.
\end{equation}
\end{lem}

\proof
Notice that
$$P_A(Y_1,\dots,Y_L) = P \left( \sum_{i=1}^L x_{i1}Y_i,\dots,\sum_{i=1}^L x_{iL}Y_i \right),$$
and so for each $v \mid \infty$,
$$H_v(P_A) \leq L^M H_v(P) \max_{1 \leq i \leq L, 1 \leq j \leq N} | x_{ij} |^{Md_v}_v \leq L^M H_v(P) \prod_{i=1}^L H_v(1,\bx_i)^M,$$
while for each $v \nmid \infty$,
$$H_v(P_A) \leq H_v(P) \max_{1 \leq i \leq L, 1 \leq j \leq N} | x_{ij} |^{Md_v}_v \leq H_v(P) \prod_{i=1}^L H_v(1,\bx_i)^M.$$
Then \eqref{ht_PA} follows by taking a product over all places of a subfield of $\qbar$ containing all the coefficients of $P$ and coordinates of $\bx_1,\dots,\bx_L$.
\endproof

Another result we require is a lemma on the existence of a small-height hyperbolic pair in a given hyperbolic plane.

\begin{lem} \label{hyper} Let $F$ be a symmetric bilinear form in $2N$ variables over $\qbar$. Let $\hyp \subseteq \qbar^N$ be a hyperbolic plane with respect to $F$. Then there exists a basis $\bx,\bwy$ for $\hyp$ such that
$$F(\bx) = F(\bwy) = 0,\ F(\bx,\bwy) \neq 0,$$
and
\begin{equation}
\label{x_bound}
H(\bx) \leq h(\bx) \leq 72\ H(F)^{\frac{1}{2}} \H(\hyp)^2,
\end{equation}
as well as
\begin{equation}
\label{y_bound}
H(\bwy) \leq h(\bwy) \leq 2592 N^2\ H(F)^{\frac{3}{2}} \H(\hyp)^4.
\end{equation}
\end{lem}

\proof
Lemma~\ref{zero} implies the existence of $\bo \neq \bx \in \hyp$ such that $F(\bx)=0$ and the height of $\bx$ is bounded as in \eqref{x_bound}. Now Theorem 1.4 of \cite{null} guarantees the existence of a point $\bz \in \hyp$ such that $F(\bz) \neq 0$ and
\begin{equation}
\label{z_bound}
H(\bz) \leq h(\bz) \leq  6 \H(\hyp).
\end{equation}
Since $F(\bx) = 0$ and $F(\bz) \neq 0$, it must be true that $\bx$ and $\bz$ are linearly independent, and hence span $\hyp$. Therefore we must have $F(\bx,\bz) \neq 0$, since $(\hyp,F)$ is regular. Then define
$$\bwy = F(\bz) \bx - 2F(\bx,\bz) \bz.$$
Clearly, $\hyp = \spn_{\qbar} \{ \bx,\bwy \}$, and it is easy to check that $F(\bwy) = 0$. Once again, regularity of $(\hyp,F)$ implies that $F(\bx,\bwy) \neq 0$, and so $\bx,\bwy$ is a hyperbolic pair basis for $\hyp$. Finally, we need to produce an estimate on the height of $\bwy$. Lemma 2.3 of \cite{smallzeros} implies that
\begin{equation}
\label{y_height}
H(\bwy) \leq h(\bwy) \leq 3N^2 H(F) h(\bx) h(\bz)^2.
\end{equation}
Combining the estimate of \eqref{y_height} with \eqref{x_bound} and \eqref{z_bound} produces \eqref{y_bound}.
\endproof

Our next lemma, which works for any field, establishes a basic divisibility property of a polynomial with respect to any fixed monomial ordering. This is Lemma~4.4 of~\cite{cfh}.

\begin{lem} \label{poly_order} Let $K$ be any field, and let $P_1(\bX), P_2(\bX) \in K[\bX] := K[X_1,\dots,X_N]$ be two polynomials in $N \geq 1$ variables over $K$. Fix any monomial ordering. Then there exist polynomials $P'_1(\bX),R(\bX) \in K[\bX]$ such that
\begin{equation}
\label{poly_lemma}
P_1(\bX) = P'_1(\bX)+R(\bX)P_2(\bX),
\end{equation}
and the leading monomial of $P_2(\bX)$, with respect to our chosen monomial ordering, does not divide any monomial of $P'_1(\bX)$.
\end{lem}

We always write $\bX$ for the variable vector $(X_1,\dots,X_N)$. Let $I \subsetneq \{1,\dots,N\}$, and write $\bX'_I$ for the vector of all variables in $\bX$ whose indices are not in $I$. The next lemma establishes the existence of zeros of especially small height for polynomials of arbitrary degree away from a hypersurface, provided the polynomial in question is of a particular form. This is an immediate adaptation of Lemma~4.5 of~\cite{cfh} over~$\qbar$ with an identical proof (word for word, while keeping in mind that the constant $A_K(d)$ of~\cite{cfh} is 1 in case of~$\qbar$), and so we do not include the proof here.

\begin{lem} \label{quad1} Let $N \geq 3$ be an integer, and let $Q(\bX) \in \qbar[\bX]$ be a polynomial of the form
\begin{equation}
\label{Q_form}
Q(\bX) = X_iX_j (c + Q_1(\bX'_{\{i,j\}})) + Q_2(\bX'_{\{i,j\}})
\end{equation}
for some indices $1 \leq i < j \leq N$, where $0 \neq c \in K$ and $Q_1$, $Q_2$ are polynomials in the $N-2$ variables $\bX'_{\{i,j\}}$. Let $P(\bX) \in \qbar[\bX]$ be a polynomial such that there exists $\bo \neq \bz \in \qbar^N$ with $Q(\bz)=0$ and $P(\bz) \neq 0$. Then there exists such $\bz$ with
\begin{equation}
\label{PQ_height}
H(\bz) \leq h(\bz) \leq H(Q).
\end{equation}
\end{lem}
\bigskip

We are now ready for the main argument of this section.

\begin{prop} \label{miss_hyper_one} Let $F$ be a quadratic form in $N$ variables over $\qbar$, and $V \subseteq \qbar^N$ be an $L$-dimensional subspace, $1 \leq L \leq N$, Suppose that the quadratic space $(V,F)$ has rank $1 \leq r \leq L$ and $\lambda$ is the dimension of the radical of $V$. Let $P(\bX) \in \qbar[X_1,\dots,X_N]$ be a polynomial of degree $M$, and suppose that there exists a nontrivial zero $\bz$ of $F$ in $V$ such that $P(\bz) \neq 0$. Then there exists such a zero $\bz$ of $F$ with
\begin{equation}
\label{z_bound_miss}
H(\bz) \leq h(\bz) \leq T(L,M) H(F)^{\max \{r,29/2\}} \H(V)^{30} 
\end{equation}
where
\begin{equation}
\label{T_const_qb}
T(L,M) = 3^{18(L-\lambda) + \frac{18L(L-1)}{L-\lambda} + \frac{33L(L-1)}{4} + 3} L^{51}.
\end{equation}
\end{prop}

\proof
First suppose that $P$ is not identically zero on $V^{\perp}$, then Theorem 1.4 of \cite{null} implies that there exists $\bo \neq \bz \in V^{\perp}$ such that $P(\bz) \neq 0$ and
$$H(\bz) \leq h(\bz) \leq 3^{\frac{\lambda(\lambda-1)}{4}} \lambda \H(V^{\perp}).$$
Combining this observation with Lemma~3.5 of~\cite{quad_qbar}, we obtain:
\begin{equation}
\label{mho_1}
H(\bz) \leq h(\bz) \leq 3^{\frac{2L(L-1) + \lambda(\lambda-1)}{4}} \lambda H(F)^r \H(V)^2,
\end{equation}
and since $F(\bz) = 0$, we are done.

Next assume that $P$ is identically zero on $V^{\perp}$. Then there must exist some nonsingular zero of $F$ on $V$ at which $P$ does not vanish; in particular, $F$ must have nonsingular zeros on $V$, so if, say, $L=1$, then we must have $V=\spn_{\qbar} \{ \bx \}$ where $F(\bx) = 0$, $P(\bx) \neq 0$ (clearly, $H(\bx) = \H(V)$ in this case), and if $L=2$, then $V$ must be a hyperbolic plane. Let $\bx_1,\dots,\bx_L$ be the small-height basis for $V$, guaranteed by Siegel's lemma (see \cite{absolute:siegel} for the original result, and Theorem~1.1 of \cite{null} for a convenient formulation):
\begin{equation}
\label{siegel_for_V}
\prod_{i=1}^L h(\bx_i) \leq 3^{\frac{L(L-1)}{4}} \H(V).
\end{equation}
Let $A = (\bx_1 \dots \bx_L)$ be the corresponding basis matrix and let $F_A$, $P_A$ be the corresponding restrictions of $F$ and $P$ to $V$ as defined in Lemma \ref{restrict}. Combining \eqref{ht_PA} with \eqref{siegel_for_V}, we obtain
\begin{equation}
\label{ht_FA}
H(F_A) \leq 3^{\frac{L(L-1)}{2}} L^2 H(F) \H(V)^2.
\end{equation}
Now notice that for each $\bz \in \qbar^L$, $F_A(\bz) = 0$, $P_A(\bz) = 0$ if and only if $F(A\bz) = 0$, $P(A\bz) = 0$, respectively. Moreover, $\bz \in \qbar^L$ is a nonsingular zero of $F_A$ if and only if $A\bz \in V$ is a non-singular zero of $F$. Also notice that by Lemma \ref{sum_height} combined with \eqref{siegel_for_V}
\begin{equation}
\label{ht_A_z}
h(A\bz) = h \left( \sum_{i=1}^L z_i \bx_i \right) \leq L h(\bz) \prod_{i=1}^L h(\bx_i) \leq 3^{\frac{L(L-1)}{4}} L h(\bz) \H(V).
\end{equation}
Since $P$ does not vanish at some nonsingular zero of $F$ on $V$, it must be that $P_A$ does not vanish at some  nonsingular zero of $F$ on $\qbar^L$; in particular, the quadratic space $(\qbar^L,F_A)$ must contain a hyperbolic plane. Our next task will be to find a hyperbolic pair of bounded height in~$(\qbar^L,F_A)$.

Lemma~3.5 of \cite{quad_qbar} states that the quadratic space $(\qbar^L, F_A)$ can be represented as $\qbar^L= \left( \qbar^L \right)^{\perp} \perp W$, where $W$ is a regular subspace of $\qbar^L$ and $\H(W) \leq 3^{\frac{L(L-1)}{2}}$. Since the quadratic spaces $(\qbar^L,F_A)$ and $(V,F)$ are isometric, their radicals have the same dimension. Therefore, the dimensions of $(\qbar^L)^{\perp}$ and $W$ are $\lambda$ and $L-\lambda$, respectively. Then Lemma~\ref{zero}  states that there exists $\bo \neq \bx \in W$ (hence $\bx$ is a nonsingular point in $(\qbar^L,F_A)$) with
\begin{eqnarray}
\label{ns_x_qb}
h(\bx) & \leq & 8 \times 3^{2(L-\lambda -1)} H(F_A)^{\frac{1}{2}} \H(W)^{\frac{4}{L-\lambda}} \nonumber \\
& \leq & 3^{2(L-\lambda) + \frac{2L(L-1)}{L-\lambda} + \frac{L(L-1)}{4}} L H(F)^{\frac{1}{2}} \H(V),
\end{eqnarray}
where the last inequality follows by \eqref{ht_FA}.

Let $\bx$ be a nonsingular point satisfying \eqref{ns_x_qb}. Since $\bx$ is nonsingular, the linear form $F_A(\bx,\bY)$ is not identically zero on $\qbar^L$, and so there must exist a standard basis vector in $\qbar^L$, call it $\bu$, such that $F_A(\bx,\bu) \neq 0$ and $h(\bu) = 1$. Then $\hyp_{xu} := \spn_\qbar \{ \bx,\bu\}$ is a hyperbolic plane in $(\qbar^L,F_A)$ with
\begin{equation}
\label{hyp_xu_nf}
\H(\hyp_{xu}) \leq L H(\bx) H(\bu) \leq 3^{2(L-\lambda) + \frac{2L(L-1)}{L-\lambda}} L^2 3^{\frac{L(L-1)}{4}} H(F)^{\frac{1}{2}} \H(V),
\end{equation}
where the first inequality is given by Lemma~\ref{lem_4.7} and the second follows by  \eqref{ns_x_qb}. Let also
$$\bwy = F_A(\bu) \bx - 2F_A(\bx,\bu) \bu,$$
then $F_A(\bwy) = 0$ and $F_A(\bx,\bwy) \neq 0$, so $\bx,\bwy$ is a hyperbolic pair. Moreover, \eqref{y_height} states that
$$h(\bwy) \leq 3L^2 H(F_A) h(\bx) h(\bu)^2 = 3L^2 H(F_A) h(\bx).$$
Combining this observation with \eqref{ns_x_qb} and \eqref{ht_FA}. we obtain that
\begin{equation}
\label{ns_y_qb}
h(\bwy) \leq 3^{2(L-\lambda) + \frac{2L(L-1)}{L-\lambda} + 1} L^5 3^{\frac{3L(L-1)}{4}} H(F)^{\frac{3}{2}}  \H(V)^3.
\end{equation}

Define
\begin{eqnarray*}
\hyp'_{xu} & := & \left\{ \bv \in \qbar^L : F_A(\bv,\bz) = 0\ \forall\ \bz \in \hyp_{xu} \right\} \\
& = & \left\{ \bv \in \qbar^L : (\bx\ \bu)^t F_A \bv = 0\ \forall\ \bz \in \hyp_{xu} \right\}
\end{eqnarray*}
to be the $(L-2)$-dimensional orthogonal complement of $\hyp_{xu}$ in $(\qbar^L,F_A)$; here we also write $F_A$ for the coefficient matrix of the quadratic form $F_A$. By the Brill-Gordan duality principle discussed in Section~\ref{notation} above, $\H(\hyp'_{xu})$ is precisely the vector space height $\H$ of the matrix $(\bx\ \bu)^t F_A$, and hence Lemma \ref{lem_2.3} implies that
$$\H(\hyp'_{xu}) \leq L^{3} H(F_A)^2 H(\bx) H(\bu),$$
and then \eqref{ht_FA} combined with \eqref{hyp_xu_nf} imply that
\begin{equation}
\label{ht_H_xu_c_2}
\H(\hyp'_{xu}) \leq 3^{2(L-\lambda) + \frac{2L(L-1)}{L-\lambda}} L^9 3^{\frac{5L(L-1)}{4}} H(F)^{\frac{5}{2}} \H(V)^5.
\end{equation}
Let $\bv_1,\dots,\bv_{L-2}$ be the small-height basis for $\hyp'_{xu}$, guaranteed by Siegel's lemma:
\begin{equation}
\label{siegel_for_hyp}
\prod_{i=1}^{L-2} h(\bv_i) \leq 3^{\frac{L^2-3L+3}{2}}\H(\hyp'_{xu}) \leq  \H(\hyp'_{xu}).
\end{equation}
Combining \eqref{siegel_for_hyp} with \eqref{ht_H_xu_c_2}, we see that
\begin{equation}
\label{ht_H_xu_c_basis_2}
\prod_{i=1}^{L-2} h(\bv_i) \leq 3^{2(L-\lambda) + \frac{2L(L-1)}{L-\lambda} + \frac{3L(L-1)}{2}} L^9 H(F)^{\frac{5}{2}} \H(V)^5.
\end{equation}
Now define the matrix $B = \left( \bx\ \bwy\ \bv_1 \dots \bv_{L-2} \right) \in \GL_L(\qbar)$, and let
$$G(\bY) = F_A(B\bY),\ Q(\bY) = P_A(B\bY).$$
Then it is easy to see that $G$ is of the form \eqref{Q_form}, and so $G$ and $Q$ satisfy the conditions of Lemma~\ref{quad1}. Hence Lemma~\ref{quad1} guarantees the existence of a point $\bw \in \qbar^L$ such that $G(\bw) = 0$, $Q(\bw) \neq 0$, and 
\begin{equation}
\label{w_point_ht}
h(\bw) \leq H(G).
\end{equation}
Now notice that standard height inequalities along with \eqref{ht_FA} imply that
\begin{eqnarray*}
H(G) & \leq & H(B^t F_A B) \leq L^{2}  H(B)^2 H(F_A) \leq L^{2} H(F_A) h(\bx)^2 h(\bwy)^2 \prod_{i=1}^{L-2} h(\bv_i)^2 \nonumber \\
& \leq & L^{4} 3^{\frac{L(L-1)}{2}} H(F) \H(V)^2 h(\bx)^2 h(\bwy)^2 \prod_{i=1}^{L-2} h(\bv_i)^2,
\end{eqnarray*}
and so by combining \eqref{ht_H_xu_c_basis_2} with \eqref{ns_x_qb} and \eqref{ns_y_qb}, we see that 
\begin{equation}
\label{G_ht_qb}
H(G) \leq 3^{12(L-\lambda) + \frac{12L(L-1)}{L-\lambda} + \frac{11L(L-1)}{2} +2} L^{34} H(F)^{10} \H(V)^{20}.
\end{equation}
Now define $\bz = A (B \bw) \in V$, and notice that $F(\bz) = F_A(B\bw) = G(\bz) = 0$, and $P(\bz) = P_A(B\bw) = Q(\bw) \neq 0$. Hence $\bz$ is precisely the point we are looking for, and to estimate its height first notice that by the same kind of reasoning as in \eqref{ht_A_z},
\begin{equation}
\label{ht_B_w_1}
h(B\bw) = h \left( w_1 \bx + w_2 \bwy + \sum_{i=1}^{L-2} w_{i+2} \bv_i \right) \leq L h(\bw) h(\bx) h(\bwy) \prod_{i=1}^L h(\bv_i).
\end{equation}
Then combining \eqref{ht_B_w_1} with \eqref{ht_A_z}, \eqref{w_point_ht}, \eqref{G_ht_qb}, \eqref{ns_x_qb}, \eqref{ns_y_qb}, and \eqref{ht_H_xu_c_basis_2} we obtain 
\begin{equation}
\label{ht_z_qb}
h(\bz) \leq T(L,M) H(F)^{10} \H(V)^{21},
\end{equation}
where $T(L,M)$ is as in \eqref{T_const_qb}. Combining \eqref{ht_z_qb} with the corresponding bound of \eqref{mho_1}, we obtain \eqref{z_bound_miss}. This completes the proof of the proposition.
\endproof

\begin{rem} \label{poly_rest} Notice that it is also easy to obtain a version of Lemma~\ref{quad1} with a restriction to a subspace $V$ of $\qbar^N$ instead of the whole $\qbar^N$ by applying Lemma \ref{restrict} in the same way as we do it in the proof of Proposition~\ref{miss_hyper_one}.
\end{rem}

\proof[Proof of Theorem \ref{miss_hyper}]
Let the notation be as in the statement of Theorem \ref{miss_hyper}. We start by extending the result of  Proposition \ref{miss_hyper_one} to a statement about a small-height zero of $F$ in $V$ outside of the union of varieties $\Z$ as defined in \eqref{Z_K}.  For our convenience, let $Z(V,F)$ be the set of nontrivial zeros of $F$ in $V$.  Since $Z(V,F) \nsubseteq \Z$, $Z(V,F) \nsubseteq Z(P_{i1},\dots,P_{ik_i})$ for all $1 \leq i \leq J$, and so for each $i$ at least one of the polynomials $P_{i1},\dots,P_{ik_i}$ is not identically zero on $Z(V,F)$, say it is $P_{ij_i}$ for some $1 \leq j_i \leq k_i$. Clearly for each $1 \leq i \leq J$, $Z(P_{i1},\dots,P_{ik_i}) \subseteq Z(P_{ij_i})$, and  $\deg(P_{ij_i}) = m_{ij_i} \leq M_i$. Define
$$P(X_1,\dots,X_N) = \prod_{i=1}^J P_{ij_i}(X_1,\dots,X_N),$$
so that $Z(V,F) \nsubseteq Z(P)$ while $\Z \subseteq Z(P)$. Then it is sufficient to construct a point of bounded height $\bx \in Z(V,F) \setminus Z(P)$. Now notice that $\deg(P) = \sum_{i=1}^J m_{ij_i} \leq M$ and apply Proposition \ref{miss_hyper_one}.

Next we want to prove the existence of a linearly independent collection of vectors $\bx_1,\dots,\bx_m \in Z(V,F) \setminus \Z$ satisfying \eqref{h_order} and \eqref{miss_hyper_bnd}, where $m=\W+\lambda$. Proposition~\ref{miss_hyper_one}, along with the argument above, guarantee the existence of a point $\bx_1 \in Z(V,F) \setminus \Z$ satisfying \eqref{z_bound_miss}. In fact, let $\bx_1$ be a point of smallest height possible in $Z(V,F) \setminus \Z$ satisfying \eqref{z_bound_miss}. If $m=1$, we are done; hence suppose that $m > 1$. Then there must exist a maximal totally isotropic subspace $W_1$ of $(V,F)$ containing $\bx_1$, and so $W_1 \nsubseteq \Z$ and $\dim_{\qbar} W_1 = m$. Then, by Theorem~A.1 of~\cite{cfh}, $W_1$ has a full basis $\bu_1,\dots,\bu_m$ outside of $\Z$. Let $\XX_1$ be an $(N-1)$-dimensional subspace of $\qbar^N$ containing $\bx_1$ so that $W_1 \nsubseteq \XX_1$, then at least one of $\bu_1,\dots,\bu_m$ is not in $\XX_1$. Since $W_1 \subseteq Z(V,F)$, we can conclude that $Z(V,F) \nsubseteq \Z^1 := \Z \cup \XX_1$, and $M(\Z^1) = M(\Z)+1$, since $\XX_1$ is the nullspace of a linear form. Again, Proposition~\ref{miss_hyper_one}, along with the argument above, guarantee the existence of a point $\bx_2 \in Z(V,F) \setminus \Z^1$ satisfying \eqref{z_bound_miss} with $M=M(\Z)+1$, and we can assume that $\bx_2$ is a point of smallest height possible in $Z(V,F) \setminus \Z^1$ satisfying \eqref{z_bound_miss}. If $m=2$, we are done; then assume $m > 2$. Then there must exist a maximal totally isotropic subspace $W_2$ of $(V,F)$ containing $\bx_1,\bx_2$, and so $W_2 \nsubseteq \Z$ and $\dim_{\qbar} W_2 = m$. Again, Theorem~A.1 of~\cite{cfh} guarantees that $W_2$ has a full basis $\bu'_1,\dots,\bu'_m$ outside of $\Z$. Let $\XX_2$ be an $(N-1)$-dimensional subspace of $V$ containing vectors $\bx_1,\bx_2$, and let $\Z^2 = \Z \cup \XX_2$. Then $V \nsubseteq \Z^2$ and $M(\Z^2) = M(\Z)+1$. Continuing to apply Proposition \ref{miss_hyper_one} and Theorem~A.1 of~\cite{cfh} in the same manner, we construct a collection of linearly independent vectors $\bx_1,\dots,\bx_m \in V \setminus \Z$ satisfying \eqref{h_order} and \eqref{miss_hyper_bnd}. This completes the proof of the theorem.
\endproof

We can now prove Corollary~\ref{miss_subspace}. First we define the constant that appears in the statement of the corollary:
\begin{eqnarray}
\label{cor_bnd_const_qb}
T^1(L,M,N,\W) & = & 3^{2(\W-1)\W^3 + \frac{(4\W+1)(L-1)(L-2)}{6}} N^{4\W+\frac{5}{2}} T(L,M+1)^{\frac{4\W+7}{3}} \nonumber \\
& \leq & 3^{2(\W-1)\W^3 + \frac{(4\W+7) L^2 (2L + 177)}{12(L-\lambda)} + (4\W+7)} N^{4\W+\frac{5}{2}} L^{17(4\W+7)}.
\end{eqnarray}

\proof[Proof of Corollary~\ref{miss_subspace}]
We first show that for each $1 \leq n \leq m=\W+\lambda$, there exists a maximal totally isotropic subspace $W^m_n$ of $(V,F)$ of bounded height, containing the corresponding point $\bx_n$ from the statement of Theorem~\ref{miss_hyper}; since $\bx_n \notin \Z$, it follows that $W^m_n \nsubseteq \Z$. First suppose that $\bx_n \in V^{\perp}$, then $\Z$ cannot contain any maximal totally isotropic subspace of $(V,F)$, since each one of them contains $V^{\perp}$. Hence we can pick $W^m_n$ to be a maximal totally isotropic subspace of $(V,F)$ of bounded height as guaranteed by Theorem~1.1 and Lemma 3.5 of \cite{quad_qbar}. Next assume that $\bx_n$ is a nonsingular point, then define
$$U_n = \{ \bz \in V : F(\bz,\bx_n) = 0 \} = \{ \bz \in K^N : \bz^t (F \bx_n) = 0 \} \cap V,$$
so that $\dim_{\qbar} U_n = L-1$ and
\begin{equation}
\label{ht_Un}
\H(U_n) \leq \H(F \bx_n) \H(V) \leq N^{3/2} H(F) H(\bx_n) \H(V),
\end{equation}
by the Brill-Gordan duality principle (discussed in Section~\ref{notation} above), combined with Lemmas~\ref{lem_2.3} and~\ref{intersection} above. Let $W'_n$ be a maximal totally isotropic subspace of $(U_n,F)$ of bounded height as guaranteed by Theorem~1.1 and Lemma 3.5 of \cite{quad_qbar}. Therefore
\begin{equation}
\label{ht_max_isot_qb}
\H(W'_n) \leq 3^{2(\W-1)\W^3 + \frac{(4\W+1)(L-1)(L-2)}{6}} H(F)^{(\W-1)^2+r} \H(U_n)^{\frac{4\W+4}{3}}.
\end{equation}
Here our bound is slightly worth than what follows from Theorem~1.1 and Lemma 3.5 of \cite{quad_qbar}, however in this form it is easier to read and apply. Now define $W^m_n=\spn_{\qbar} \{ \bx_n,W'_n \}$, then $W^m_n$ is a maximal totally isotropic subspace of $(V,F)$ containing $\bx_n$. Moreover,
$$\H(W^m_n) \leq N^{1/2}\ H(\bx_n) \H(W'_n).$$
Combining this observation with \eqref{ht_max_isot_qb}, and the bounds of Theorem~\ref{miss_hyper}, we obtain \eqref{cor_bnd} in the case $k=m$.

Siegel's lemma implies the existence of a basis $\bw_1,\dots,\bw_m$ for $W^m_n$ such that
$$\prod_{i=1}^m h(\bw_i) \leq 3^{\frac{m(m-1)}{4}} \H(W^m_n).$$
Since $\bo \neq \bx_n \in W^m_n$, there must exist a subcollection of $m-1$ of these vectors which are linearly independent with $\bx_n$; since we did not order these vectors by height, we can assume without loss of generality that $\bx_n,\bw_2,\dots,\bw_m$ are linearly independent. Then for each $1 \leq k < m$, define
$$W^k_n = \spn_{\qbar} \{ \bx_n,\bw_2,\dots,\bw_k \},$$
so that $\bx_n \in W^k_n$, $\dim_{\qbar} W^k_n = k$,
$$\spn_{\qbar} \{ \bx_n \} = W^1_n \subset W^2_n \subset \dots \subset W^m_n,$$
and by Lemma \ref{lem_4.7}
$$\H(W^k_n) \leq N^{k/2} H(\bx_n) \prod_{i=2}^k h(\bw_i) \leq 3^{\frac{m(m-1)}{4}} N^{k/2} H(\bx_n) \H(W^m_n),$$
which is precisely \eqref{cor_bnd_1}. This completes the proof of the corollary.
\endproof
\bigskip

\section{Simultaneous zeros of a system of quadratic forms}
\label{section_k_forms}

The main goal of this section is to prove Theorem~\ref{k_forms}. We first prove a more technical version of this result, from which the theorem is then derived.

\begin{prop} \label{k-forms} Let $k \geq 1$ be an integer, $F_1,\dots,F_k$ be quadratic forms in $N$ variables over $\qbar$, and let $V \subseteq \qbar^N$ be an $L$-dimensional subspace, $N \geq L \geq \frac{k(k+1)}{2} + 1$. There exists $\bo \neq \bw \in V$ such that $F_m(\bw) = 0$ for all $1 \leq m \leq k$ and
\begin{equation}
\label{hw}
h(\bw) \leq \left( 3^{\frac{L^2}{2}} N^{\frac{3(L+1)}{2}} \H(V) \right)^{a_k} \prod_{m=1}^k H(F_m)^{b_k(m)},
\end{equation}
where $b_1(1)=1/2$, $b_k(k)=2$ for all $k \geq 2$, and
\begin{equation}
\label{bkm}
b_k(m) = 2\left( (k+2) b_{k-1}(m) + 1 \right) + 4 \left( (k-1) b_{k-1}(m) + 1 \right) \sum_{m=1}^{k-1} b_{k-1}(m),
\end{equation}
as well as $a_1 = 2$, and
\begin{equation}
\label{ak}
a_k = 2\left( (k+1) a_{k-1} + 1 \right) + 4 \left( (k-1) a_{k-1} + 1 \right) \sum_{m=1}^{k-1} b_{k-1}(m).
\end{equation}
\end{prop}

\proof
We argue by induction on $k$. If $k=1$, the result is given by Lemma~\ref{zero}. By induction hypothesis, we can pick $\bo \neq \bx \in V$ such that $F_m(\bx) = 0$ for all $1 \leq m \leq k-1$ and
\begin{equation}
\label{bxx}
h(\bx) \leq \left( 3^{\frac{L^2}{2}} N^{\frac{3(L+1)}{2}} \H(V) \right)^{a_{k-1}} \prod_{m=1}^{k-1} H(F_m)^{b_{k-1}(m)}.
\end{equation}
If $F_k(\bx) = 0$, we are done, so assume $F_k(\bx) \neq 0$. Let
$$W = \left\{ \bwy \in V : F_m(\bx,\bwy) = 0\ \forall\ 1 \leq m \leq k-1 \right\},$$
then $\dim_{\qbar} W \geq L - (k-1) \geq \frac{k(k-1)}{2} + 2$ and
\begin{equation}
\label{WFm}
\H(W) \leq  N^{\frac{3(k-1)}{2}} H(\bx)^{k-1} \H(V) \prod_{m=1}^{k-1} H(F_m),
\end{equation}
by Lemmas~\ref{lem_2.3} and~\ref{intersection}. Let $\ell = \dim_{\qbar} W - 1 \geq \frac{k(k-1)}{2} + 1$. Let $\bwy_1,\dots,\bwy_{\ell+1}$ be the small-height basis for $W$, guaranteed by Siegel's lemma, then
\begin{equation}
\label{yi}
\prod_{i=1}^{\ell+1} h(\bwy_i) \leq 3^{\frac{L(L-1)}{2}} \H(W) \leq 3^{\frac{L(L-1)}{2}} N^{\frac{3(k-1)}{2}} H(\bx)^{k-1} \H(V) \prod_{m=1}^{k-1} H(F_m),
\end{equation}
where the second inequality follows from~\eqref{WFm}. Notice that $\bx \in W$, then at least $\ell$ vectors of the $\bwy_i$'s above are linearly independent with~$\bx$; let $\bwy_1,\dots,\bwy_\ell$ be these vectors. Let
\begin{equation}
\label{w_ab}
\bw = \alpha \bx + \sum_{i=1}^{\ell} \beta_i \bwy_i \in V,
\end{equation}
where values $\alpha,\beta_1,\dots,\beta_{\ell} \in \qbar$, not all zero, are to be specified. The first observation is that $\bw \neq 0$. For each $1 \leq m \leq k-1$, define
$$g_m(\beta_1,\dots,\beta_{\ell}) := F_m(\bw) = \sum_{i=1}^{\ell} \sum_{j=1}^{\ell} F_m(\bwy_i,\bwy_j) \beta_i \beta_j,$$
since $F_m(\bx,\bwy_i) = 0$ for all $1 \leq i,j \leq \ell$. By Lemma~\ref{F_xy}, for each place $v \in M(\qbar)$,
\begin{eqnarray*}
H_v(g_m) & = & \max_{1 \leq i,j \leq \ell} |F_m(\bwy_i,\bwy_j)|^{d_v}_v \\
& \leq & \left\{ \begin{array}{ll}
N^{2d_v} H_v(F_m) \max_{1 \leq i,j \leq \ell} H_v(\bwy_i) H_v(\bwy_j) & \mbox{if $v \mid \infty$,} \\
H_v(F_m) \max_{1 \leq i,j \leq \ell} H_v(\bwy_i) H_v(\bwy_j) & \mbox{if $v \nmid \infty$.}
\end{array}
\right.
\end{eqnarray*}
Then
\begin{eqnarray}
\label{gm_bnd}
H(g_m) & \leq & N^2 H(F_m) \prod_{i=1}^{\ell} h(\bwy_i)^2 \nonumber \\
& \leq & 3^{L(L-1)} N^{3k-1} H(\bx)^{2(k-1)} \H(V)^2 H(F_m) \prod_{n=1}^{k-1} H(F_n)^2,
\end{eqnarray}
where the last inequality follows by~\eqref{yi}. Notice that $\H(\qbar^{\ell})=1$. Then, by induction hypothesis, there exists $\bo \neq \bb = (b_1,\dots,b_{\ell}) \in \qbar^{\ell}$ such that
$$g_1(\bb) = g_2(\bb) = \dots = g_{k-1}(\bb)$$
and
\begin{eqnarray}
\label{hbb}
h(\bb) & \leq & \prod_{m=1}^{k-1} H(g_m)^{b_{k-1}(m)} \nonumber \\
& \leq & \left( 3^{L(L-1)} N^{3k-1} H(\bx)^{2(k-1)} \H(V)^2 \prod_{m=1}^{k-1} H(F_m)^2\right)^{\sum_{m=1}^{k-1} b_{k-1}(m)} \nonumber \\
& \times & \prod_{m=1}^{k-1} H(F_m)^{b_{k-1}(m)}
\end{eqnarray}
Set $\beta_i=b_i$ for each $1 \leq i \leq \ell$ in~\eqref{w_ab}, then
\begin{equation}
\label{F_k}
F_k(\bw) = \alpha^2 F_k(\bx) + 2\alpha \sum_{i=1}^{\ell} F_k(\bx,\bwy_i) b_i + \sum_{i=1}^{\ell} \sum_{j=1}^{\ell} F_k(\bwy_i,\bwy_j) b_i b_j.
\end{equation}
Setting \eqref{F_k} equal 0 and solving for $\alpha$, we obtain 
\begin{eqnarray*}
\alpha & = & \frac{2}{F_k(\bx)} \Bigg( - \sum_{i=1}^{\ell} F_k(\bx,\bwy_i) b_i  \\
& \ & \pm \left\{ \left( \sum_{i=1}^{\ell} F_k(\bx,\bwy_i) b_i \right)^2 - F_k(\bx)  \sum_{i=1}^{\ell} \sum_{j=1}^{\ell} F_k(\bwy_i,\bwy_j) b_i b_j \right\}^{\frac{1}{2}} \Bigg).
\end{eqnarray*}
We now need to estimate the Weil height of $\alpha$. First notice that, similar to Lemma~\ref{F_xy},
\begin{equation}
\label{F_xx}
\max \left\{ 1, \left| \frac{2}{F_k(\bx)} \right|^{d_v}_v \right\} \leq  \left\{ \begin{array}{ll}
(2N^2)^{d_v} H_v(1, F_k) H_v(1,\bx)^2  & \mbox{if $v \mid \infty$,} \\
H_v(1,F_k) H_v(1,\bx)^2 & \mbox{if $v \nmid \infty$.}
\end{array}
\right.
\end{equation}
Suppose $v \mid \infty$, then applying triangle inequality we obtain:
\begin{eqnarray}
\label{sum_ineq-1}
& & \left| - \sum_{i=1}^{\ell} F_k(\bx,\bwy_i) b_i \pm \left\{ \left( \sum_{i=1}^{\ell} F_k(\bx,\bwy_i) b_i \right)^2 - F_k(\bx)  \sum_{i=1}^{\ell} \sum_{j=1}^{\ell} F_k(\bwy_i,\bwy_j) b_i b_j \right\}^{\frac{1}{2}} \right|_v  \\
& & \leq \ell \max_{1 \leq i \leq \ell} |b_i|_v \left( 2 \max_{1 \leq i \leq \ell} |F_k(\bx,\bwy_i)|_v + |F_k(\bx)|_v^{\frac{1}{2}} \max_{1 \leq i,j \leq \ell} |F_k(\bwy_i,\bwy_j)|^{\frac{1}{2}}_v \right) \nonumber \\
& & \leq 3 \ell N^2 \max_{1 \leq i \leq \ell} |b_i|_v \left( H_v(F_k) H_v(\bx) \max_{1 \leq i \leq \ell} H_v(\bwy_i) \right)^{\frac{1}{d_v}}, \nonumber
\end{eqnarray}
where the last inequality follows by Lemma~\ref{F_xy}. Similarly, applying the ultrametric inequality along with Lemma~\ref{F_xy} in case $v \nmid \infty$, we obtain:
\begin{eqnarray}
\label{sum_ineq-2}
& & \left| - \sum_{i=1}^{\ell} F_k(\bx,\bwy_i) b_i \pm \left\{ \left( \sum_{i=1}^{\ell} F_k(\bx,\bwy_i) b_i \right)^2 - F_k(\bx)  \sum_{i=1}^{\ell} \sum_{j=1}^{\ell} F_k(\bwy_i,\bwy_j) b_i b_j \right\}^{\frac{1}{2}} \right|_v \\
& & \leq \max_{1 \leq i \leq \ell} |b_i|_v \left( H_v(F_k) H_v(\bx) \max_{1 \leq i \leq \ell} H_v(\bwy_i) \right)^{\frac{1}{d_v}}. \nonumber
\end{eqnarray}
Finally, notice that dividing through by one of the nonzero coefficients, if necessary, we can assume without loss of generality that
$$H_v(F_k) = H_v(1,F_k)$$
for every $v \in M(K)$. Combining \eqref{F_xx}, \eqref{sum_ineq-1}, \eqref{sum_ineq-2}, and taking a product over all $v \in M(\qbar)$, we have:
\begin{equation}
\label{a_bnd}
h(\alpha) \leq 6 \ell N^4 H(F_k)^2 h(\bx)^3 h(\bb) \prod_{i=1}^{\ell} h(\bwy_i).
\end{equation}
Then, applying Lemma~\ref{sum_height} to $\bw$ and using inequalities \eqref{a_bnd} and \eqref{yi}, along with observation that $\ell \leq L-1$, we obtain:
\begin{eqnarray}
\label{w-bnd}
h(\bw) & \leq & (\ell+1) h(\alpha) h(\bb) h(\bx) \prod_{i=1}^{\ell} h(\bwy_i) \nonumber \\
& \leq & 6 L^2 3^{L(L-1)} N^{3k+1} h(\bx)^{2(k+1)} h(\bb)^2 \H(V)^2 \prod_{m=1}^k H(F_m)^2.
\end{eqnarray}
Now~\eqref{hw} follows by combining~\eqref{w-bnd} with~\eqref{bxx} and~\eqref{hbb}.
\endproof

We can now establish the bound of Theorem~\ref{k_forms}.

\proof[Proof of Theorem~\ref{k_forms}] Let the notation be as in Proposition~\ref{k-forms}. Define
$$B_k = \max_{1 \leq m \leq k-1} b_k(m),$$
then~\eqref{bkm} implies that
\begin{equation}
\label{B_k}
B_k = 2\left( (k+2) B_{k-1} + 1 \right) + 4 \left( (k-1) B_{k-1} + 1 \right) (k-1) B_{k-1},
\end{equation}
where $B_1=b_1(1)=1/2$, hence $B_2=9$. Then~\eqref{B_k} implies that for $k \geq 3$,
$$B_k \leq \left( 2k B_{k-1} \right)^2,$$
and so
\begin{equation}
\label{BK-1}
B_k \leq \frac{1}{4} \times 36^{2^{k-2}} \prod_{m=3}^k m^{2^{k-m+1}}.
\end{equation}
Also notice that, by~\eqref{ak},
$$a_k = 2\left( (k+1) a_{k-1} + 1 \right) + 4 \left( (k-1) a_{k-1} + 1 \right) (k-1) B_{k-1},$$
where $a_1=2$, hence $a_2=20$. Then for $k \geq 3$,
$$a_k \leq 4k^2 a_{k-1} B_{k-1} \leq a_{k-1} 36^{2^{k-2}} k^2 \prod_{m=3}^k m^{2^{k-m+1}}$$
where the last inequality follows by~\eqref{BK-1}. Therefore:
$$a_k \leq \frac{20}{36^2} \times 36^{2^{k-1}} \left( \prod_{m=3}^k m^2 \right) \left( \prod_{m=3}^k m^{2(2^{k-m+1}-1)} \right) = \frac{20}{81} B_k^2.$$
This completes the proof.
\endproof
\bigskip

\section{Zeros of inhomogeneous quadratic polynomials}
\label{inh}

In this section we consider the inhomogeneous situation. In particular, we obtain bounds for the height of zeros of a system of one or two inhomogeneous quadratic polynomials and a collection of inhomogeneous linear equations over~$\qbar$, thus proving Theorem~\ref{thm_inh}.

\begin{prop} \label{inh-1} Let $F$ be a quadratic polynomial in $N \geq 2$ variables over $\qbar$, possibly inhomogeneous. Let $m$ be an integer, $0 \leq m \leq N-2$, and $\L_1,\dots,\L_m$ be linear polynomials in $N$ variables over $\qbar$, possibly inhomogeneous; the case $m=0$ just means that there are no linear polynomials. Suppose that the system
\begin{equation}
\label{system}
F(\bx) = \L_1(\bx) = \dots = \L_m(\bx) = 0
\end{equation}
has a nontrivial solution over~$\qbar$. Then there exists such a solution $\bo \neq \bwy \in \qbar^N$ with
\begin{equation}
\label{hw_mn}
h(\bwy) \leq 8 (N+1)^{2m} 3^{2(N-m+1)(N-m)} H(F)^{\frac{1}{2}} \prod_{i=1}^m H(\L_i)^4.
\end{equation}
\end{prop}

\proof
Introduce one more variable $x_{N+1}$ to homogenize the system~\eqref{system}. These new forms in $N+1$ variables, which we will still denote by $F$ and $\L_j$, $1 \leq i \leq k$, $1 \leq j \leq m$ (if $m > 0$), have the same heights as the corresponding original polynomials. Let
\begin{equation}
\label{V_def}
V = \left\{ \bx \in \qbar^{N+1} : \L_i(\bx) = 0\ \forall\ 0 \leq i \leq m \right\},
\end{equation}
which is a subspace of $\qbar^{N+1}$ of dimension $N+1-m$ and, by Lemma~\ref{intersection}, the Brill-Gordan duality principle, and ~\eqref{ht_ineq_sqrt},
\begin{equation}
\label{H-V}
\H(V) \leq \H(\L_1) \cdots \H(\L_m) \leq (N+1)^{m/2} \prod_{i=1}^m H(\L_i).
\end{equation}
Then $\bx \in \qbar^N$ is a solution of~\eqref{system} if and only if there exists some $x_{N+1} \neq 0$ such that $\bx' := (\bx,x_{N+1}) \in V$ is a zero of the form $F'$. Clearly, $h(\bx) \leq h(\bx')$. Define 
$$U = \left\{ \bx \in \qbar^{N+1} : x_{N+1} = 0 \right\},$$
then $\H(U) = 1$.

Our argument here is an adaptation of the proof of Lemma 4.1 of \cite{quad_qbar} to the situation of an inhomogeneous quadratic polynomial. Let $W_1 = V \cap U$. Since a solution to~\eqref{system} exists, it must be true that 
$$2 \leq \dim W_1 = N-m < \dim V = N+1-m.$$
Let $\bwy_1,\dots,\bwy_{N+1-m}$ be the small-height basis for $V$ guaranteed by Siegel's lemma. At least one of these vectors must be in $V \setminus W_1$, call it $\bz_1$. Then
\begin{equation}
\label{z1}
h(\bz_1) \leq 3^{\frac{(N-m+1)(N-m)}{2}}  \H(V) \leq 3^{\frac{(N-m+1)(N-m)}{2}} (N+1)^{\frac{m}{2}} \prod_{i=1}^m H(\L_i).
\end{equation}
If $F(\bz_1) = 0$, we are done, so assume $F(\bz_1) \neq 0$. For each vector $\bv \in V$, there must exist $\bu \in W_1$ and $\alpha \in \qbar$ such that
$$\bv = \alpha \bz_1 + \bu.$$
Since there exists some $\bv \in V$ like this for which
$$F(\bv) = \alpha^2 F(\bz_1) + 2\alpha F(\bz_1,\bu) + F(\bu) = 0,$$
it must be true that either $F(\bz_1,\bu) \neq 0$ or $F(\bu) \neq 0$ for some $\bu \in W_1$. Now Theorem 1.4 of \cite{null} guarantees the existence of a point $\bz_2 \in W_1$ such that either $F(\bz_1,\bz_2) \neq 0$ or $F(\bz_2) \neq 0$ and
\begin{equation}
\label{z2_bound}
h(\bz_2) \leq  2 \times 3^{\frac{(N-m+1)(N-m)}{2}} \H(W_1),
\end{equation}
where $\H(W_1)=\H(V)$ by Lemma~\ref{intersection}. Therefore
\begin{equation}
\label{z12}
h(\bz_1) h(\bz_2) \leq 2 (N+1)^{m} 3^{(N-m+1)(N-m)} \prod_{i=1}^m H(\L_i)^2.
\end{equation}
We will now construct $0 \neq a_1,a_2 \in \qbar$ such that $\bwy = a_1 \bz_1 + a_2 \bz_2 \in V$ is a small-height zero of $F$. By construction of $\bz_1,\bz_2$, we will have $y_{N+1} \neq 0$. We want
\begin{equation}
\label{newform}
0 = F(\bwy) = F(\bz_1) a_1^2 + 2F(\bz_1,\bz_2) a_1 a_2 + F(\bz_2) a_2^2 = G(a_1,a_2).
\end{equation}
The right hand side of (\ref{newform}) is a quadratic form $G$ in the variables $a_1,a_2$ with coefficients $F(\bz_1), 2F(\bz_1,\bz_2),F(\bz_2)$. Notice that either $(a_1,a_2) = (0,0)$ or both $a_1,a_2 \neq 0$. By Lemma~3.4 of~\cite{quad_qbar}, there must exist such a pair $(a_1,a_2) \neq (0,0)$ with
\begin{equation}
\label{q0}
h(a_1,a_2) \leq 2 \sqrt{H(G)}.
\end{equation}
Let $E$ be the field extension generated over $K$ by coefficients of $G$. By Lemma~\ref{F_xy}, for each $v \in M(E)$, we have
\begin{eqnarray}
\label{q1}
H_v(G) & \leq & \max \{ |F(\bz_1)|_v, |2|_v |F(\bz_1,\bz_2)|_v, |F(\bz_2)|_v \} \nonumber \\
& \leq & H_v(F) \max \{ H_v(\bz_1)^2, H_v(\bz_1)H_v(\bz_2), H_v(\bz_2)^2 \} \\
& \leq & H_v(F) \max \{ 1, H_v(\bz_1) \}^2 \max \{ 1, H_v(\bz_2) \}^2. \nonumber
\end{eqnarray}
if $v \nmid \infty$, and
\begin{eqnarray}
\label{q2}
H_v(G)^{\frac{2d}{d_v}} & \leq & \|F(\bz_1)\|^2_v + 2 \|F(\bz_1,\bz_2)\|^2_v + \|F(\bz_2)\|^2_v \leq H_v(F)^{\frac{2d}{d_v}} \times \nonumber \\
& \times & \left( H_v(\bz_1)^{\frac{4d}{d_v}} + 2 \left( H_v(\bz_1) H_v(\bz_2) \right)^{\frac{2d}{d_v}} + H_v(\bz_2)^{\frac{4d}{d_v}} \right) \\
& \leq & H_v(F)^{\frac{2d}{d_v}} \left( 1 + H_v(\bz_1)^{\frac{2d}{d_v}} \right)^2 \left( 1 + H_v(\bz_2)^{\frac{2d}{d_v}} \right)^2. \nonumber
\end{eqnarray}
if $v | \infty$.
Combining (\ref{q0}) with (\ref{q1}) and (\ref{q2}), we see that
\begin{equation}
\label{q3}
h(a_1,a_2) \leq 2 H(F)^{\frac{1}{2}} h(\bz_1) h(\bz_2).
\end{equation}
Combining \eqref{q3} and Lemma~\ref{sum_height}, we see that there exists a zero of $F$ of the form $\bwy = a_1 \bz_1 + a_2 \bz_2 \in V$ so that 
\begin{equation}
\label{q4}
h(\bwy) \leq 2 H(F)^{\frac{1}{2}} h(\bz_1)^2 h(\bz_2)^2.
\end{equation}
Now \eqref{hw_mn} follows by combining \eqref{z12} and \eqref{q4}. This completes the proof.
\endproof

\begin{rem} \label{better} The result of Proposition~\ref{inh-1} also readily follows from our Theorem~\ref{miss_hyper} as a small special case, but with a weaker bound.
\end{rem}

\begin{prop} \label{inh-2} Let $F$ and $G$ be quadratic polynomials in $N \geq 4$ variables over $\qbar$, possibly inhomogeneous. Let $m$ be an integer, $0 \leq m \leq N-4$, and $\L_1,\dots,\L_m$ be linear polynomials in $N$ variables over $\qbar$, possibly inhomogeneous; the case $m=0$ just means that there are no linear polynomials. Suppose that the system
\begin{equation}
\label{system-1}
F(\bx) = G(\bx) = \L_1(\bx) = \dots = \L_m(\bx) = 0
\end{equation}
has a nontrivial solution over~$\qbar$. Then there exists such a solution $\bo \neq \bw \in \qbar^N$ with
\begin{equation}
\label{hw_mn-1}
h(\bw) \leq \M(m,N) H(F)^{58} H(G)^3 \prod_{i=1}^m H(\L_i)^{180},
\end{equation}
where
\begin{equation}
\label{MmN}
\M(m,N) = 18 \times 8^{38} (N+1)^{90m+8} (N+1-m)^{36} 3^{90(N-m+1)(N-m)}.
\end{equation}
\end{prop}

\proof
We start as in the proof of Proposition~\ref{inh-1} above, homogenizing our system~\eqref{system-1} with an additional variable $x_{N+1}$ and defining $V$ as in~\eqref{V_def}. Assuming the existence of a simultaneous zero $\bwy \in V$ of $F$ and $G$ with $y_{N+1} \neq 0$, it is our goal now to produce such a point $\bwy$ of bounded height.

Let $\bx \in V$ be such that $x_{N+1} \neq 0$ and $F(\bx) = 0$, satisfying~\eqref{hw_mn}, as guaranteed by Proposition~\ref{inh-1}. If $G(\bx) = 0$, we are done, so assume $G(\bx) \neq 0$. Define
$$W = \left\{ \bwy \in V : F(\bx,\bwy) = 0 \right\},$$
then 
$$4 = (N+1)-(N-4)-1 \leq \dim V - 1 \leq \dim W \leq \dim V = N+1-m.$$
Let
$$W_1 = \left\{ \bu \in W : u_{N+1} =0 \right\},\ W_2 = \spn_{\qbar} \{ \bx \}.$$
Clearly, $W_1,W_2 \subsetneq W$, and so Theorem~A.1 of \cite{cfh} guarantees the existence of linearly independent vectors $\bwy,\bz \in W \setminus (W_1 \cup W_2)$ such that
\begin{eqnarray}
\label{y-z-ht}
\max \{ h(\bwy), h(\bz) \} & \leq & 3^{\frac{(N-m)(N+1-m)}{2}}  (N+1-m)  \H(W) \nonumber \\ 
& \leq & 3^{\frac{(N-m)(N+1-m)}{2}} (N+1-m) (N+1)^{\frac{m}{2}} h(\bx) H(F) \prod_{i=1}^m H(\L_i),
\end{eqnarray}
where the last inequality follows by Lemmas~\ref{intersection} and~\ref{lem_2.3} along with~\eqref{H-V}. We will now construct a common zero of $F$ and $G$ with a nonzero last coordinate of the form
$$\bw = a \bx + b \bwy + c \bz,$$
where the choice of $a,b,c$ is specified below, and so
\begin{equation}
\label{wxyz}
h(\bw) \leq 3 h(a,b,c) h(\bx) h(\bwy) h(\bz),
\end{equation}
by Lemma~\ref{sum_height}. Notice that $F(\bw) = F(b \bwy + c \bz)$, and set $F(\bw) = G(\bw) = 0$ and $w_{N+1}=1$, then we obtain a system of three equations in the coefficients $a,b,c$:
\begin{equation}
\label{FG_abc}
\left. \begin{array}{ll}
b^2 F(\bwy) + 2bc F(\bwy,\bz) + c^2 F(\bz) = 0 \\
a^2 G(\bx) + b^2 G(\bwy)  + c^2 G(\bz) + 2(ab G(\bx,\bwy) + ac G(\bx,\bz) + bc G(\bwy,\bz)) = 0 \\
a x_{N+1} + b y_{N+1} + c z_{N+1} = 1
\end{array}
\right\}.
\end{equation}
\smallskip

{\it Case 1.} Suppose first that every triple $a,b,c$ satisfying the second equation of the system~\eqref{FG_abc} has $a=0$. This implies that $G(\bwy) = G(\bz) = G(\bwy,\bz) = 0$, and hence $G$ is identically zero on $U_1 := \spn_{\qbar} \{ \bwy,\bz \}$. If $F(\bwy) = 0$ or $F(\bz) = 0$, take $\bw = \bwy$ or $\bz$, respectively, and the result follows by combining~\eqref{y-z-ht} with Proposition~\ref{inh-1}. Hence assume $F(\bwy) F(\bz) \neq 0$. Taking
$$b = \frac{1-cz_{N+1}}{y_{N+1}},$$
we can set
$$\bw = \frac{1-cz_{N+1}}{y_{N+1}} \bwy + c \bz,$$
and so $w_{N+1} \neq 0$ for any choice of $c$. To choose $c$, let
\begin{eqnarray}
\label{Fwyz}
F(\bw) & = & \left( \frac{1-cz_{N+1}}{y_{N+1}} \right)^2 F(\bwy) + \left( \frac{2c(1-cz_{N+1})}{y_{N+1}} \right) F(\bwy,\bz) + c^2 F(\bz) \nonumber \\
& = & \frac{1}{y_{N+1}^2} \left\{ F(z_{N+1} \bwy - y_{N+1} \bz) c^2 - 2 F(\bwy, z_{N+1} \bwy - y_{N+1} \bz) c + F(\bwy) \right\} = 0.
\end{eqnarray}
Since $F(\bwy) \neq 0$, it must be true that at least one of  $F(z_{N+1} \bwy - y_{N+1} \bz)$ and $F(\bwy, z_{N+1} \bwy - y_{N+1} \bz)$ is not equal to zero. If $F(z_{N+1} \bwy - y_{N+1} \bz) = 0$, take
\begin{equation}
\label{c_id1}
c = \frac{F(\bwy)}{2F(\bwy, z_{N+1} \bwy - y_{N+1} \bz)}.
\end{equation}
If $F(z_{N+1} \bwy - y_{N+1} \bz) \neq 0$, take
\begin{eqnarray}
\label{c_id2}
c & = & \frac{2}{F(z_{N+1} \bwy - y_{N+1} \bz)} \Big( F(\bwy, z_{N+1} \bwy - y_{N+1} \bz) \nonumber \\
& & \pm \sqrt{ F(\bwy, z_{N+1} \bwy - y_{N+1} \bz)^2 - F(z_{N+1} \bwy - y_{N+1} \bz) F(\bwy)} \Big).
\end{eqnarray}
Applying Lemmas~\ref{sum_height}, \ref{F_xy} and taking the maximum of heights of the right hand sides of~\eqref{c_id1} and~\eqref{c_id2}, we see that
$$h(c) \leq 24 (N+1)^4 H(F)^2 h(\bwy)^4 h(\bz)^4 h(z_{N+1},y_{N+1})^4.$$
Therefore
\begin{eqnarray}
\label{ht_w_1}
h(\bw) & \leq & 2 h \left( c, \frac{1-cz_{N+1}}{y_{N+1}} \right) h(\bwy) h(\bz) \leq 2 h(c)^2 h(y_{N+1},z_{N+1}) h(\bwy) h(\bz) \nonumber \\
& \leq & 1152 (N+1)^8 H(F)^4 h(\bwy)^{18} h(\bz)^{18} \nonumber \\
& \leq & \M_1(m,N) H(F)^{40} h(\bx)^{36} \prod_{i=1}^m H(\L_i)^{36},
\end{eqnarray}
where the last inequality follows by~\eqref{y-z-ht} and the constant
$$\M_1(m,N) = 1152 (N+1)^{18m+8} (N+1-m)^{36} 3^{18(N-m)(N+1-m)}.$$
Combining~\eqref{ht_w_1} with Proposition~\ref{inh-1}, we obtain
\begin{equation}
\label{ht_w_2}
h(\bw) \leq 8^{36} \M_1(m,N) (N+1)^{72m} 3^{72(N-m+1)(N-m)} H(F)^{58}  \prod_{i=1}^m H(\L_i)^{180}.
\end{equation}
\smallskip

{\it Case 2.} Now suppose that there exists some triple $(a,b,c)$ with $a \neq 0$ satisfying the equations~\eqref{FG_abc} and assume that $F(\bwy)=0$ (respectively, $F(\bz)=0$). Then $F$ is identically zero on $U_2 := \spn_{\qbar} \{ \bx,\bwy \}$ (respectively, on $U_3 := \spn_{\qbar} \{ \bx,\bz \}$, and we can repeat the argument from Case 1 above for the form $G$ on $U_2$ (respectively, on $U_3$) instead of $F$ on $U_1$. The bound we obtain on the height of the resulting point $\bw$ is smaller than that of~\eqref{ht_w_2}.
\smallskip

{\it Case 3.} Next suppose that there exists some triple $(a,b,c)$ with $a \neq 0$ satisfying the equations~\eqref{FG_abc} and $F(\bwy)F(\bz) \neq 0$. Then we can solve this system for $a,b,c$, and estimating the height of such solution, we obtain:
$$h(a,b,c) \leq 2016 (N+1)^{16} H(F)^5 H(G)^3 h(\bx)^4 h(\bwy)^4 h(\bz)^4.$$
Applying Lemma~\ref{sum_height}, we see that
\begin{eqnarray}
\label{w_a_n0}
h(\bw) & \leq & 6048 (N+1)^{16} H(F)^5 H(G)^3 h(\bx)^5 h(\bwy)^5 h(\bz)^5 \nonumber \\
& \leq & \M_2(m,N) H(F)^{15} H(G)^3 h(\bx)^{15} \prod_{i=1}^m H(\L_i)^{10} \nonumber \\
& \leq & 8^{15} (N+1)^{30m} 3^{30(N-m+1)(N-m)} \M_2(m,N) H(F)^{\frac{45}{2}} H(G)^3 \prod_{i=1}^m H(\L_i)^{70},
\end{eqnarray}
where
$$\M_2(m,N) = 6048 (N+1)^{5m+16} (N+1-m)^{10} 3^{5(N-m)(N+1-m)},$$
and the last two inequalities follow by~\eqref{y-z-ht} and Proposition~\ref{inh-1}, respectively. Inequality~\eqref{hw_mn-1} is now obtained by combining~\eqref{ht_w_2} with~\eqref{w_a_n0}. This finishes the proof of the proposition.
\endproof

\proof[Proof of Theorem~\ref{thm_inh}] The result of the theorem now follows by combining Propositions~\ref{inh-1} and~\ref{inh-2}.
\endproof
\bigskip

{\bf Acknowledgment.} The author thanks Wai Kiu Chan and Glenn R. Henshaw for many discussions and their helpful comments on the subject of this paper, and acknowledges the wonderful hospitality and support of the Erwin Schr\"odinger Institute for Mathematical Physics in Vienna, Austria, where a part of this work was done. The author is also grateful to the referee for the helpful remarks and suggestions.
\bigskip

\bibliographystyle{plain}  
\bibliography{quad_zero-qbar}        

\end{document}